\documentclass{amsart}

\usepackage{graphicx} 

\newtheorem{LM}{Lemma}[section]
\newtheorem{THM}[LM]{Theorem}
\newtheorem{Q}[LM]{Question}

\newtheorem{CL}[LM]{Claim}
\newtheorem{CON}[LM]{Conjecture}

\begin{document}

\title[Non primitive/Seifert-fibered construction]
{Seifert fibered surgeries which do not arise from 
primitive/Seifert-fibered constructions} 

\author{Thomas Mattman$^*$, \ 
Katura Miyazaki 
and Kimihiko Motegi$^{**}$ \\
}
\address{Department of Mathematics and Statistics, 
California State University, Chico, Chico CA95929-0525, U.S.A. }
\email{TMattman@CSUChico.edu}
\address{Faculty of Engineering, Tokyo Denki University, Tokyo 101-8457, 
Japan}
\email{miyazaki@cck.dendai.ac.jp}
\address{Department of Mathematics, Nihon University, 
Tokyo 156-8550, Japan}
\email{motegi@math.chs.nihon-u.ac.jp}


\date{}

\maketitle

\begin{abstract}
We construct two infinite families of knots each of which admits
a Seifert fibered surgery 
with none of these surgeries coming 
from Dean's primitive/Seifert-fibered construction.
This disproves a conjecture that all Seifert fibered 
surgeries arise 
from Dean's primitive/Seifert-fibered construction. 
The $(-3,3,5)$-pretzel knot belongs to both of 
the infinite families.
\end{abstract}

{%
\renewcommand{\thefootnote}{}
\footnotetext{$^{*}$Supported in part by grants
from NSERC and FCAR}
\footnotetext{$^{**}$Supported in part by Grant-in-Aid for 
Scientific Research (No.\ 40219978),
The Ministry of Education, Culture, Sports, Science and Technology, Japan.}
}%

\section{Introduction}

        Let $K$ be a knot in the $3$-sphere $S^3$.
Then we denote by $(K;\gamma)$ the $3$-manifold obtained by 
$\gamma$-surgery on $K$, 
i.e., by attaching a solid torus to $S^3-$int$N(K)$ in such a way 
that $\gamma$ bounds a meridian disk of the filling solid torus. 
Using the preferred meridian-longitude pair of $K\subset S^3$,
we parametrize slopes $\gamma$ of $K$
by $r \in \mathbb{Q} \cup \{ \infty \}$;
then we also write $(K; r)$ for $(K; \gamma)$. \par

We begin by recalling Berge's \cite{B} construction, 
an explicit
construction which yields several infinite families of knots each
admitting a lens space Dehn surgery. \par

        Let $K$ be a knot contained in a genus two Heegaard surface 
$F$ for $S^3$, 
i.e., 
$S^3 = H \cup_{F} H'$, 
where $H$ and $H'$ denote genus two handlebodies. 
Suppose that $K$ is nontrivial and 
that the manifolds $H(K)$ and $H'(K)$ are both solid tori,
where $H(K)$ (resp. $H'(K)$) is obtained by
attaching  a $2$-handle to $H$ (resp. $H'$) along $K$. 
The isotopy class in $\partial N(K)$ of the curve(s) in 
$\partial N(K) \cap F$ 
is called the \textit{surface slope} of $K$ 
with respect to $F$. 
Then by performing Dehn surgery on $K$ along 
the surface slope $\gamma$, 
we obtain a $3$-manifold 
$(K; \gamma) = H(K) \cup H'(K)$,
which is a lens space.
It cannot be $S^2 \times S^1$ by \cite{Ga},
nor $S^3$ by \cite{GL}. 
This construction is called 
\textit{Berge's construction} or the
\textit{primitive/primitive construction} 
and such a knot $K$ is said to be \textit{primitive/primitive} 
with respect to $F$. \par

        In \cite{B} Berge suggested the following. 
See also \cite{Go}. 

\begin{CON}
\label{Berge}
If $(K; \gamma)$ is a lens space, 
then this surgery arises from Berge's construction. 
\end{CON}

Dean \cite{D}, \cite{D2} made a natural modification to Berge's
construction; 
suppose that $K$ is as before except that $H'(K)$ is now 
a Seifert fiber space over the disk with two exceptional fibers. 
Then for the surface slope $\gamma$, 
$(K; \gamma)$ is a Seifert fiber space over 
$S^2$ with at most three exceptional fibers
or a connected sum of two lens spaces. 
If $K$ is hyperbolic, 
then the cabling conjecture \cite{GS} 
states that the latter cannot occur. 
This construction is called 
\textit{Dean's construction}
or the \textit{primitive/Seifert-fibered construction} 
and such a knot $K$ is said to be 
\textit{primitive/Seifert-fibered} 
with respect to $F$. \par

        The notion of primitive/Seifert-fibered construction 
has been slightly generalized by allowing the possibility
that $H'(K)$ is a Seifert fiber space 
over the M{\"o}bius band 
with one exceptional fiber \cite{EM2}, \cite{MM7}. 
In the following, we use the term 
primitive/Seifert-fibered construction (or knot) 
in this generalized sense. \par

        In analogy with Conjecture \ref{Berge}, 
Dean \cite{D} and Gordon \cite{Go} asked:

\begin{Q}
\label{Dean}
If $(K; \gamma)$ is a Seifert fiber space other than a lens space, 
then does this surgery arise 
from a primitive/Seifert-fibered 
construction?  
\end{Q}

        Many examples of Seifert fibered surgeries 
(see, for example, \cite{BH}, \cite{BZ}, \cite{EM1} and \cite{EM2}) 
have been constructed using the Montesinos trick 
(\cite{Mon}, \cite{Bl}). 
        Recently, in \cite{EM2}, 
Eudave-Mu\~noz has shown that all known examples 
of Seifert fibered surgeries constructed by the Montesinos trick 
can be explained by Dean's construction. 
Furthermore, Seifert fibered surgeries on 
twisted torus knots in \cite{MM3} can also be explained by 
such constructions \cite{MM7}. 

        On the other hand, 
in the present note we demonstrate the following which 
answers the question above in the negative. 
A knot $K$ is \textit{strongly invertible} if there is an orientation 
preserving involution of $S^3$ which leaves $K$ invariant and reverses 
an orientation of $K$; 
primitive/Seifert-fibered knots are shown to be strongly invertible. \par

\begin{THM}
\label{non-p/s}
        There is an infinite family of non-strongly invertible knots
each of which 
admits a Seifert fibered surgery with none of these surgeries 
arising from the primitive/Seifert-fibered construction. 
For example, the $(-3,3,5)$-pretzel knot belongs to the family.
\end{THM}

	Very recently Hyung-Jong Song has observed that the 1-surgery of the
$(-3,3,3)$-pretzel knot is a Seifert fibered surgery, 
but does not arise from the primitive/Seifert-fibered construction.
In contrast with our examples, the $(-3,3,3)$-pretzel knot is
strongly invertible; 
but it has cyclic period $2$ and 
tunnel number greater than one like ours.

In his thesis \cite{Mat},
the first author observed that the $(-3, 3, 5)$-pretzel knot 
has a small Seifert fibered surgery by 
experiments via Weeks' computer program SnapPea.  
This observation is the starting point of our study. 

\medskip 

\textbf{Acknowledgements} --
The first author wishes to thank Steven Boyer and Jinha Jun 
for helpful conversations. 
We would like to thank the referee for careful reading and useful comments. 

\section{Examples}

        We shall say that a Seifert fiber space is of 
\textit{type} $S^2(n_1, n_2, n_3)$ 
if it has a Seifert fibration over 
$S^2$ with three exceptional fibers of 
indices $n_1, n_2$ and $n_3$ $(n_i \ge 2)$.  

\bigskip

\begin{figure}[htbp]
\begin{center}
\includegraphics[width=0.3\linewidth]{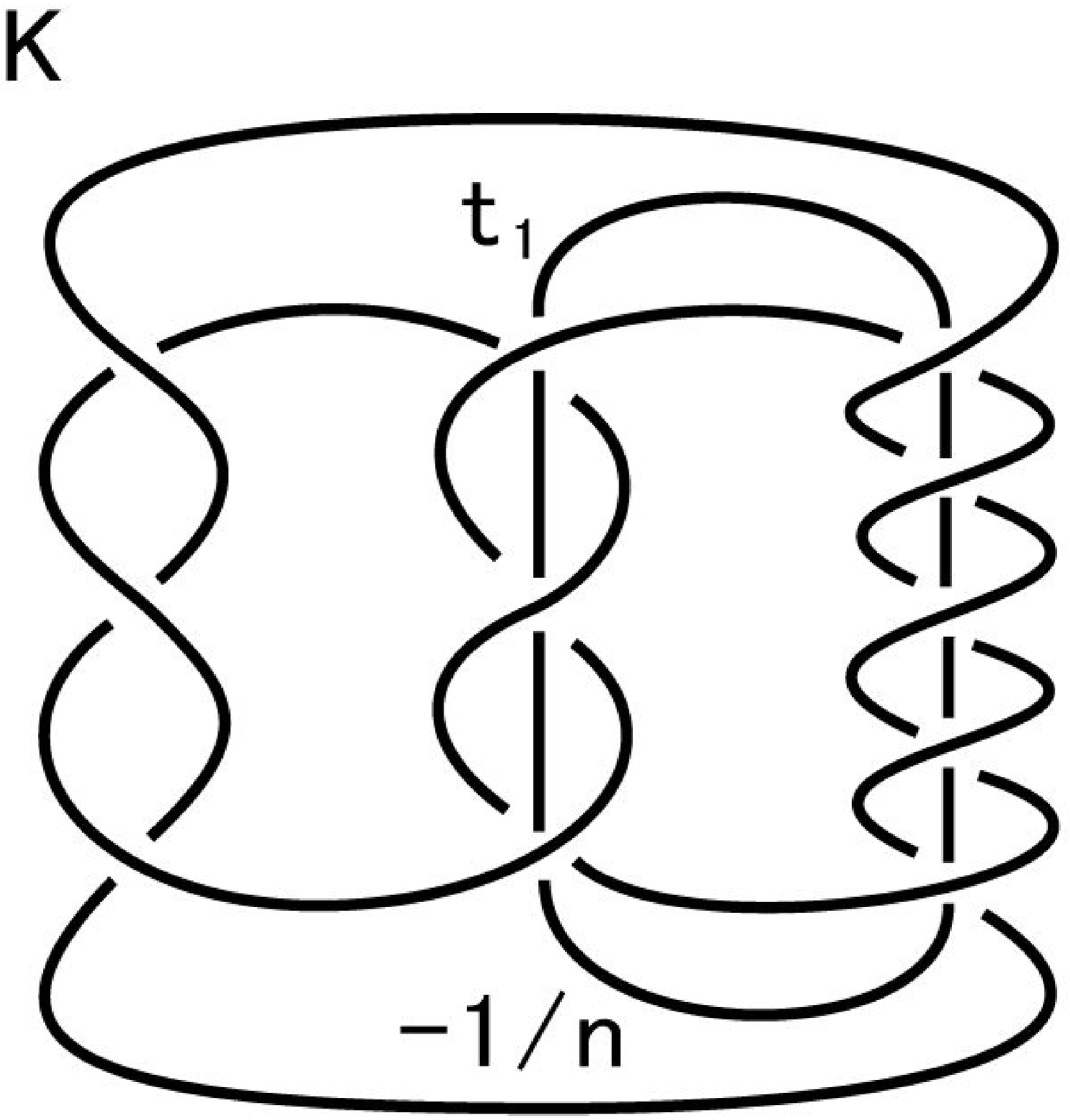}
\caption{}
\label{fig:Fig1}
\end{center}
\end{figure}

\noindent
\textbf{Example 1.}
        Let $K \cup t_1$ be the two component link of Figure \ref{fig:Fig1}. 

Here $K$ is the Montesinos knot given by the triple of 
rational tangles $(1/3, -1/3, -1/5)$, 
which is often called
the $(-3, 3, 5)$-pretzel knot. 
(We adopt Bleiler's convention \cite{Ble} 
on the parametrization of rational tangles.) 
Let $K_n$ ($n$ is possibly zero) be the knot obtained from $K$ by performing 
$-1/n$-surgery on $t_1$. 
Equivalently, 
$K_n$ is obtained by doing $n$-twisting along $t_1$. 
Then $K_n$ enjoys the following properties. 
\begin{enumerate}
\item $K_n$ is a hyperbolic knot, 
\item $K_n$ has cyclic period $2$, 
but is not strongly invertible, 
\item the tunnel number of $K_n$ is $2$, and 
\item $(K_n ; 1)$ is a Seifert fiber space of type 
$S^2(3, 5, |15n + 4|)$. 
\end{enumerate}

\medskip

Before verifying properties (1)--(4) 
we observe that $\{K_n\}$ is the family of Theorem~\ref{non-p/s}.

\medskip

\begin{figure}[htbp]
\begin{center}
\includegraphics[width=0.9\linewidth]{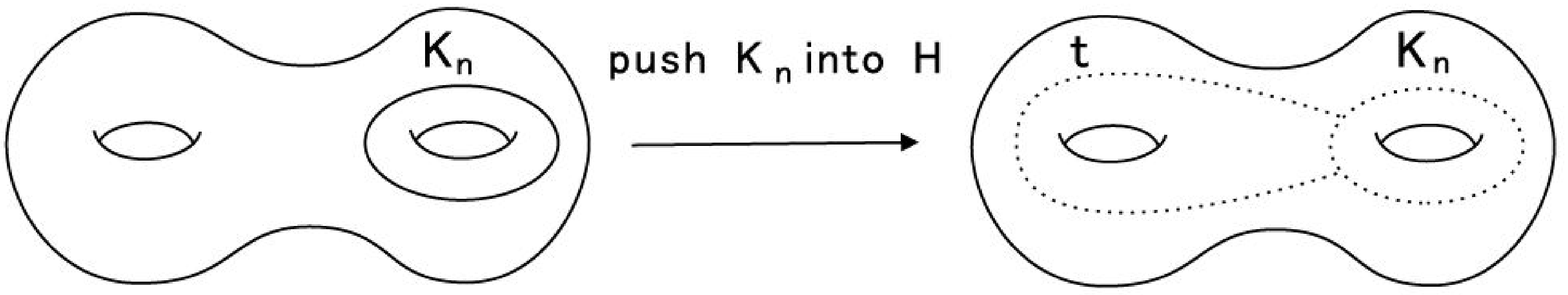}
\caption{}
\label{fig:Fig2}
\end{center}
\end{figure}

\noindent
\textit{Proof of Theorem} \ref{non-p/s}.
Properties (2) and (4) show that $K_n$ is not strongly invertible and
admits a Seifert fibered surgery. 
Assume for a contradiction that $K_n$ is primitive/Seifert-fibered; 
then $H(K_n)$ is a solid torus for an unknotted genus $2$ handlebody 
$H$ with $K \subset \partial H$. 
First we show that $K_n$ has tunnel number $1$ following \cite{D}. 
By \cite{Z}, there is a homeomorphism of the genus two handlebody $H$ 
after which $K_n$ appears as in Figure \ref{fig:Fig2}.
After pushing $K_n$ into $H$, 
take an arc $t$ as in Figure \ref{fig:Fig2}. 

\noindent
Then $H - \mathrm{int}N(K_n \cup t)$ is the product of a surface 
and an interval. 
Thus $S^3 - \mathrm{int}N(K_n \cup t) = 
H' \cup (H - \mathrm{int}N(K_n \cup t))$ is a genus two handlebody, 
so the knot $K_n$ has tunnel number $1$. 
This then implies that $K_n$ is strongly invertible by 
\cite[Lemma 5]{Mor2}, a contradiction. 
Hence the Seifert fibered surgery does not come from the 
primitive/Seifert-fibered construction.
\hspace*{\fill} $\qed$(Theorem~\ref{non-p/s})

\begin{CL}
\label{Period2}
$K_n$ has cyclic period $2$.
\end{CL}

\begin{figure}[htbp]
\begin{center}
\includegraphics[width=0.4\linewidth]{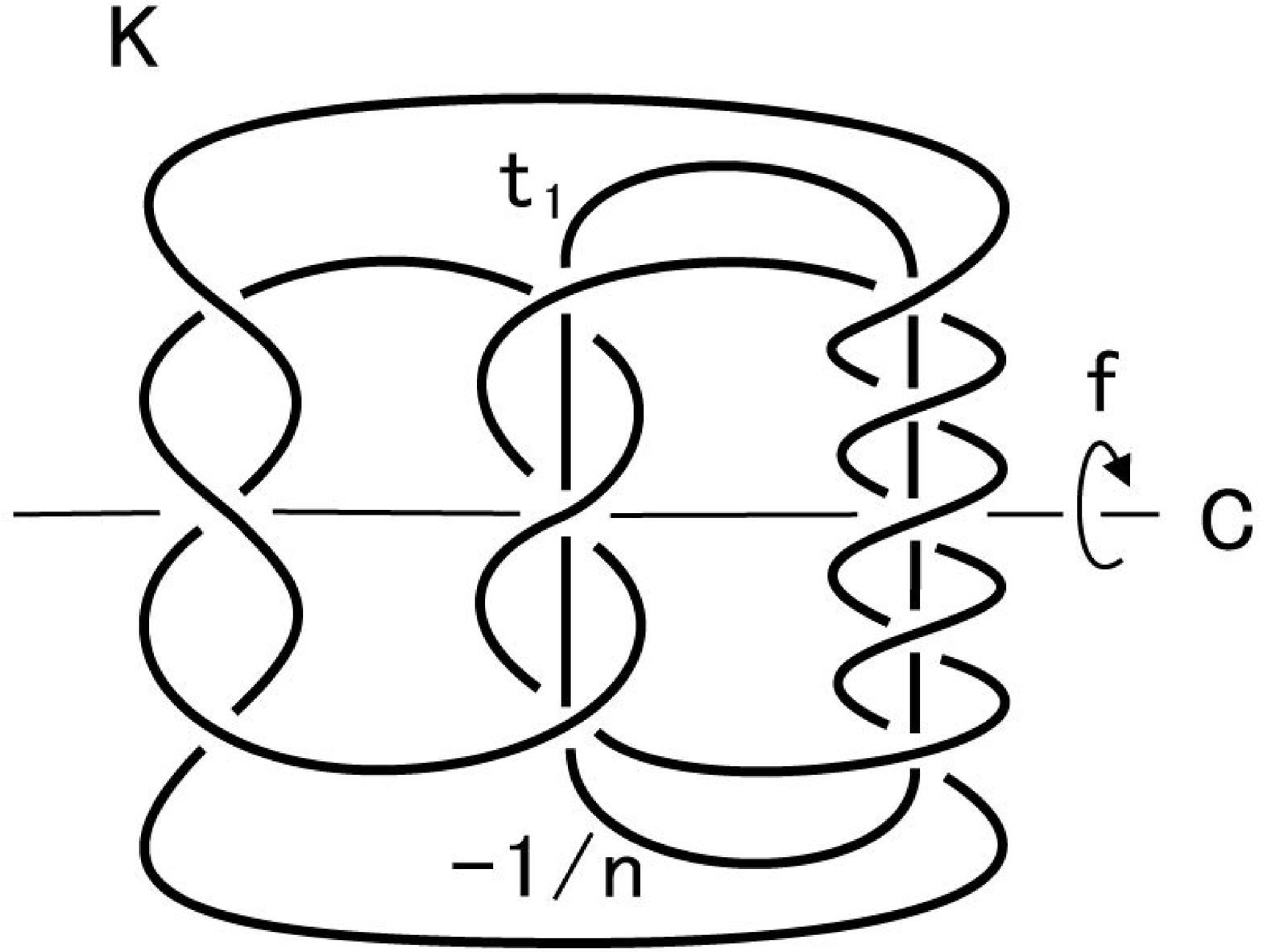}
\caption{}
\label{fig:Fig3}
\end{center}
\end{figure}

\noindent 
\textit{Proof.}
	As shown in Figure \ref{fig:Fig3}, 
let $f: S^3 \to S^3$ be
the $\pi$-rotation about $C$ such that $f(K) = K$ and $f(t_1) = t_1$.
The axis $C$ is disjoint from $K$
and intersects $t_1$ in exactly two points. 

Hence, $f|S^3-$int$N(t_1)$ extends to an involution $\bar{f}$
of $\displaystyle (t_1; -1/n)\cong S^3$ about an axis $\overline{C}$
such that $\bar{f}(K_n) =K_n$ and $K_n\cap \overline{C} =\emptyset$.
It follows that $K_n$ has cyclic period 2.
\hspace*{\fill} $\qed$(Claim~\ref{Period2})

\begin{CL}
\label{Seifert}
$(K_n; 1)$ is a Seifert fiber space of type 
$S^2(3, 5, |15n + 4|)$.
\end{CL}

\noindent
\textit{Proof.}
        Let $(K \cup t_1;\ 1, -1/n)$ denote the manifold 
obtained by performing a surgery on the link 
$K \cup t_1$ with surgery slopes $1$ for $K$ and 
$-1/n$ for $t_1$. 
We will show that $(K \cup t_1;\ 1, -1/n)$ is a 
Seifert fiber space of type 
$S^2(3, 5, |15n + 4|)$. \par

\begin{figure}[htbp]
\begin{center}
\includegraphics[width=0.75\linewidth]{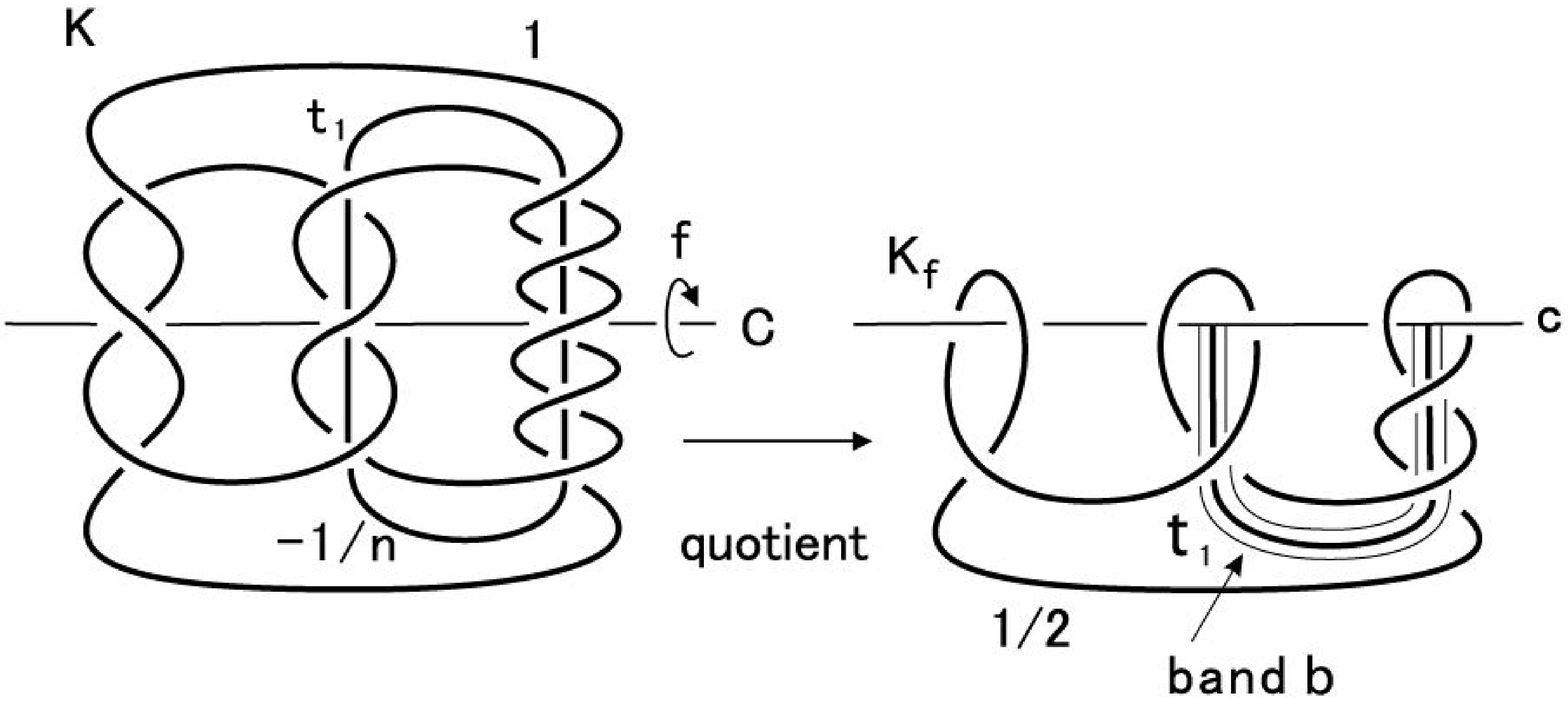}
\caption{}
\label{fig:Fig4}
\end{center}
\end{figure}

        To prove this 
we form the quotient by the involution $f:S^3 \to S^3$ to 
obtain the factor knot $K_f$, 
the branched knot $c$ which is the image of $C$, 
and the arc $\tau_1$ which is the 
image of $t_1$ and connects 
two points in $c$ (Figure \ref{fig:Fig4}). 

As shown in Figure \ref{fig:Fig4}, 
the factor knot $K_f$ is unknotted in $S^3/f \cong S^3$. 
Note that $1$-surgery on $K$ corresponds to $1/2$-surgery on 
the factor knot $K_f$ which is equivalent to 
$(-2)$-twisting along $K_f$ because $K_f$ is unknotted;
see Figure \ref{fig:Fig6}.   
We denote the image of $c$ after $(-2)$-twisting along $K_f$ by $c'$. 
Note also that by the Montesinos trick (\cite{Mon}, \cite{Bl}), 
$-1/n$-surgery on $t_1$ corresponds to 
$-1/n$-untangle surgery 
(i.e., a replacement of a $1/0$-untangle by 
a $-1/n$-untangle) on $c'$ along $\tau_1$ 
as indicated in Figure \ref{fig:Fig8}. 
In order to correctly perform the untangle surgery, 
we keep track of the framing. 
This can be done by indicating a band $\beta$ whose core 
is $\tau_1$; see Figure \ref{fig:Fig4}.  
(For simplicity, 
we indicate the band $\beta$ in only two places: 
just after taking 
the quotient by the involution $f$, and 
just before performing the untangle surgery.)
By an isotopy as in Figures \ref{fig:Fig6} and \ref{fig:Fig7}, 
we see that $c'$ is the Montesinos knot 
given by the triple of rational tangles 
$(2/5, -3/4, 1/3)$. 
Denote the result of $-1/n$-untangle surgery on 
$c'$ by $c'_n$ (Figure \ref{fig:Fig8}). 
Then $c'_n$ is the Montesinos knot 
given by the triple of rational tangles
$(2/5, (11n+3)/(-15n-4), 1/3)$, 
and the branched covering space 
$(K \cup t_1;\ 1, -1/n)$ of $S^3$ 
branched along $c'_n$ 
is a Seifert fiber space of type 
$S^2(3, 5, |15n+4|)$. 
Since the linking number of $K$ and $t_1$ is zero, 
the $1$-slope of $K$ corresponds to the $1$-slope of $K_n$, 
and hence $(K \cup t_1;\ 1, -1/n) \cong (K_n; 1)$. 
It follows that $(K_n; 1)$ is a Seifert fiber space
of type $S^2(3, 5, |15n+4|)$ as required.  
\hspace*{\fill} $\qed$(Claim \ref{Seifert})

\begin{figure}[htbp]
\begin{center}
\includegraphics[width=1.0\linewidth]{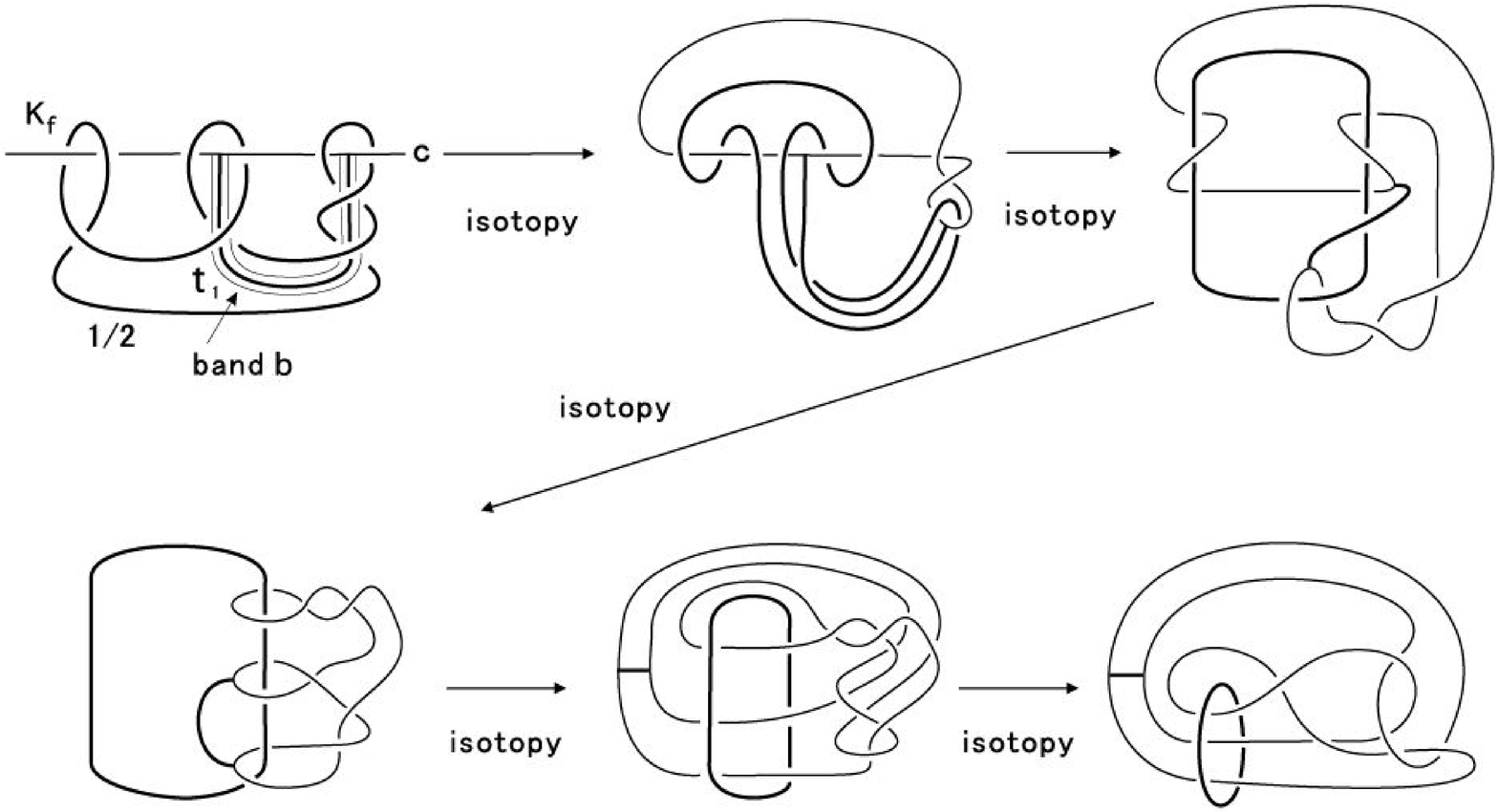}
\caption{}
\label{fig:Fig5}
\end{center}
\end{figure}

\begin{figure}[htbp]
\begin{center}
\includegraphics[width=1.0\linewidth]{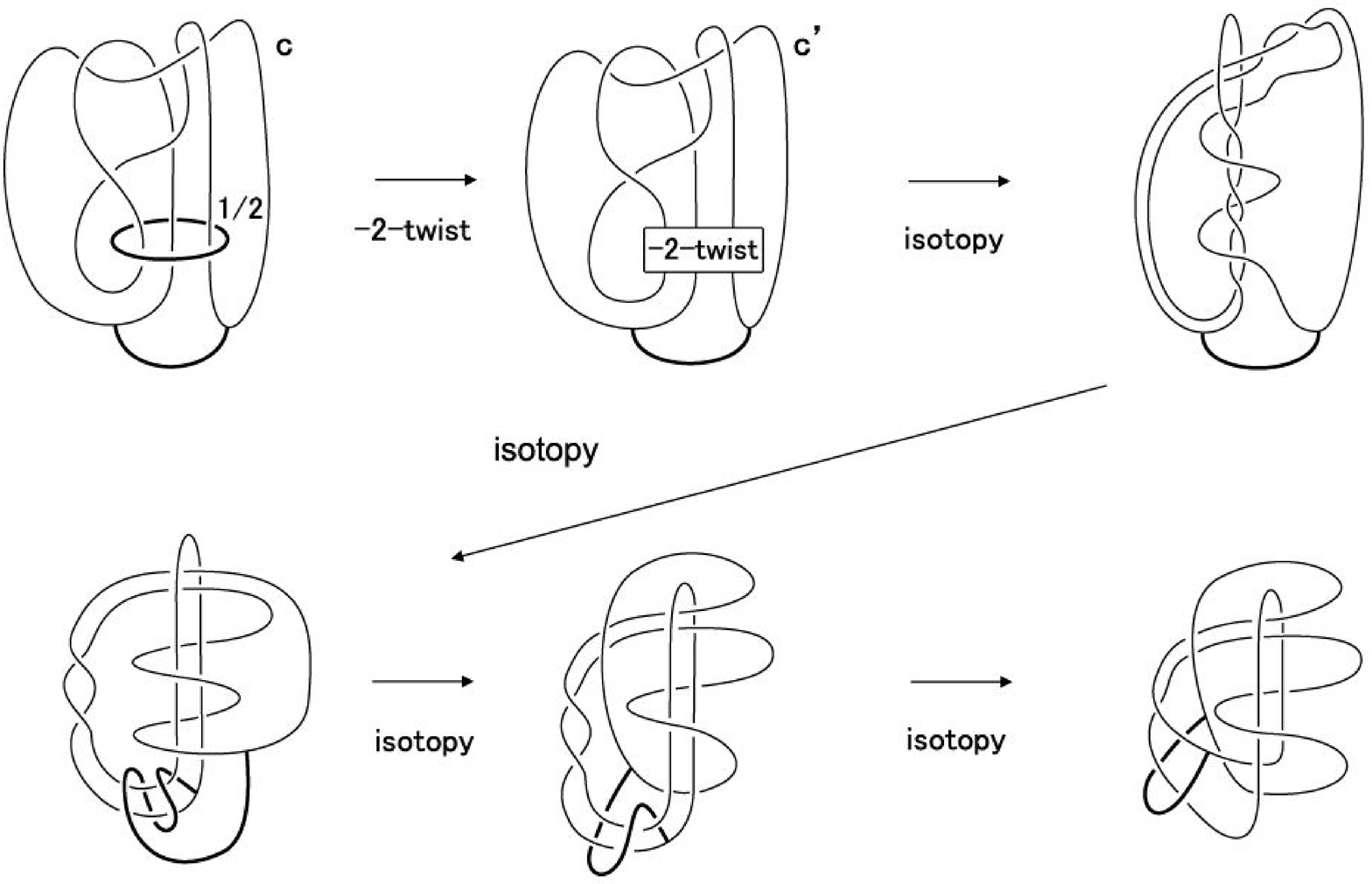}
\caption{Continued from Figure \ref{fig:Fig5}.}
\label{fig:Fig6}
\end{center}
\end{figure}

\begin{figure}[htbp]
\begin{center}
\includegraphics[width=1.0\linewidth]{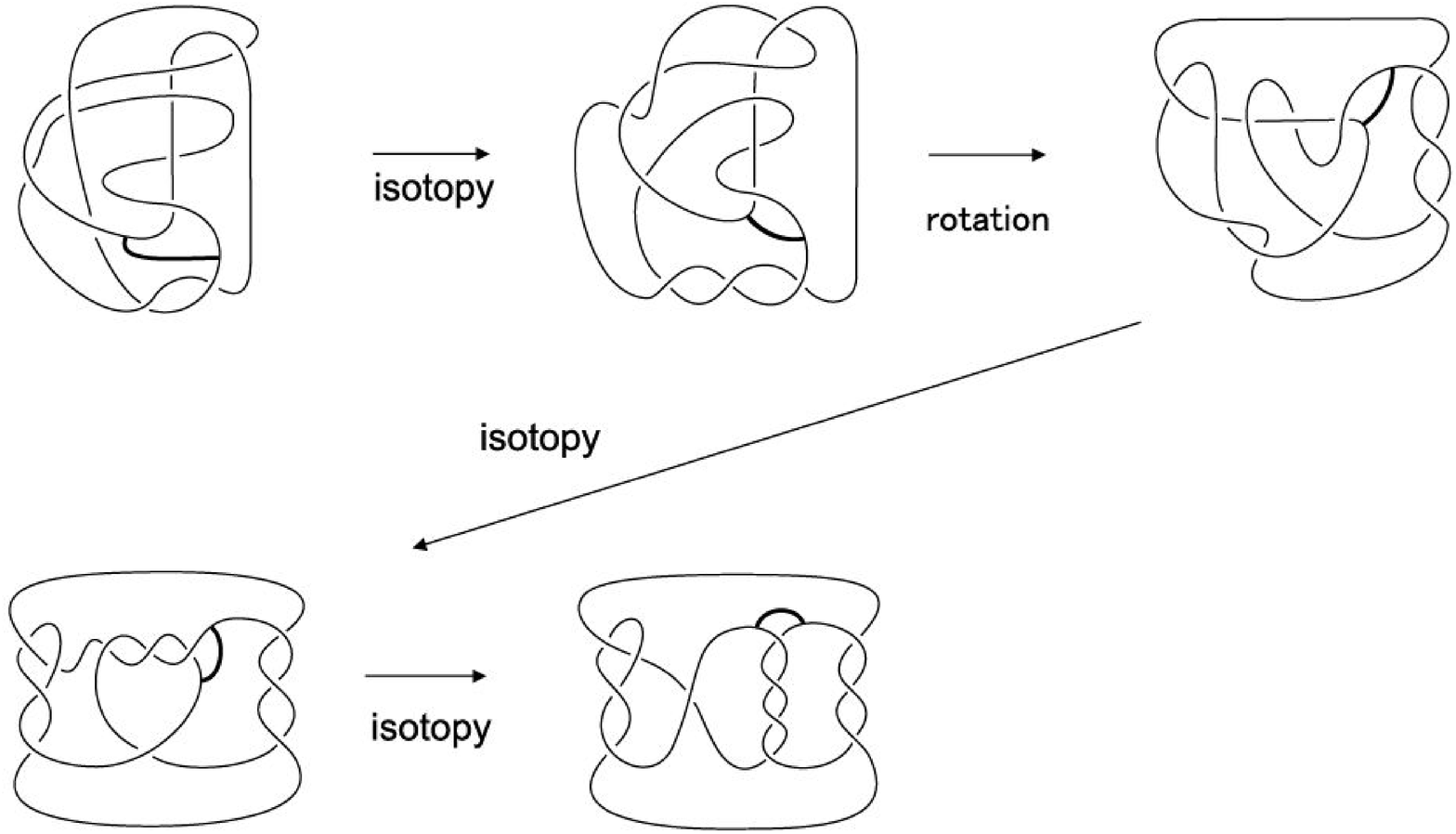}
\caption{Continued from Figure \ref{fig:Fig6}.}
\label{fig:Fig7}
\end{center}
\end{figure}

\begin{figure}[htbp]
\begin{center}
\includegraphics[width=1.0\linewidth]{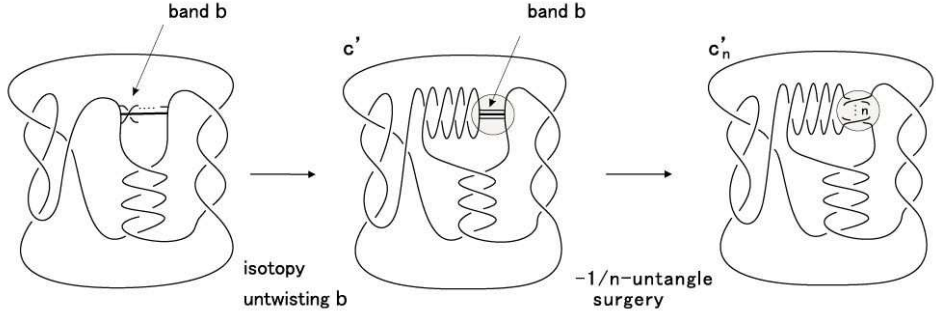}
\caption{Continued from Figure \ref{fig:Fig7}.}
\label{fig:Fig8}
\end{center}
\end{figure}

\medskip

\begin{CL}
\label{hyperbolic}
$K_n$ is a hyperbolic knot. 
\end{CL}

\noindent
\textit{Proof.}
The knot $K$ bounds an obvious Seifert surface $S$ of genus one.
Since $t_1$ can be isotoped off $S$, after doing $n$-twisting
along $t_1$ $S$ becomes a Seifert surface for $K_n$.
By Claim \ref{Seifert}, 
$K_n$ is a nontrivial knot and 
thus $g(K_n)$, 
the genus of $K_n$, 
is equal to one. \par

	Assume for a contradiction that $K_n$ is a satellite knot. 
Then since $(K_n ; 1)$ is atoroidal, 
$K_n$ has a companion solid torus $V$ whose core is a simple knot 
$\widehat{K_n}$ such that $K_n$ is a 
$0$ or $1$-bridge braid in $V$ (\cite[Proposition 2.2(1)]{MM2}). 
From Schubert's formula \cite{S} 
(\cite[Proposition 2.10]{BuZ}) we have 
$g(K_n) \ge w g(\widehat{K_n})$, 
where $w$ denotes the winding number of $K_n$ in $V$. 
Since $w \ge 2$ and $g(\widehat{K_n}) \ge 1$, 
we have $g(K_n) \ge 2$, a contradiction. 
        If $K_n$ is a torus knot, 
then since the genus is one, 
$K_n$ is a $(\pm 2, 3)$-torus knot $T_{\pm 2, 3}$. 
However $(T_{2, 3}; 1)$ (resp. $(T_{-2, 3}; 1)$) is 
a Seifert fiber space of type 
$S^2(2, 3, 5)$ (resp. $S^2(2, 3, 7)$), 
contradicting Claim \ref{Seifert}. 
It follows that $K_n$ is a hyperbolic knot. 
\hspace*{\fill} $\qed$(Claim \ref{hyperbolic})

\begin{CL}
\label{non-invertible}
$K_n$ is not strongly invertible. 
\end{CL}

\noindent
\textit{Proof.}
	Recall that $K_n$ has cyclic period $2$ and that 
$(K_n; 1)$ is a Seifert fiber space of type 
$S^2(3, 5, |15n+4|)$ (Claim \ref{Seifert}).
Since $|15n + 4| > 2$ and  $|15n + 4|  \ne 3, \ 5$, 
if $K_n$ is strongly invertible, 
then by \cite[Theorem 1.7(1)]{Mot}, 
$K_n$ is a torus knot or a cable of a torus knot. 
This contradicts $K_n$ being hyperbolic (Claim \ref{hyperbolic}). 
Therefore $K_n$ is not strongly invertible. 
\hspace*{\fill} $\qed$(Claim \ref{non-invertible})

\begin{CL}
\label{tunnel number}
The tunnel number of $K_n$ is two. 
\end{CL}

\noindent
\textit{Proof}.
Let $H$ be a handlebody in $S^3$ 
which is obtained by thickening the obvious genus one Seifert surface
for $K$.
Then $F =\partial H$ is a genus 2 Heegaard surface for $S^3$
which contains $K$.
Since $t_1$ is a core of a handlebody $H$,
$H$ remains a handlebody after $-1/n$-surgery on $t_1$.
It follows that $K_n$ is embedded in a genus 2 Heegaard surface $F$.
Then, by \cite[Fact on p.138]{Mor}
the tunnel number of $K_n$ is less than or equal to $2$.
On the other hand, 
since a tunnel number one knot is strongly invertible
(\cite[Lemma 5]{Mor2}),  
Claim \ref{non-invertible} implies that 
the tunnel number of $K_n$ is two. 
\hspace*{\fill} $\qed$(Claim \ref{tunnel number})

\bigskip

\noindent
\textbf{Example 2.} 
The second example is a variant of Example 1. 
Let us consider the trivial knot $t_2$ 
of Figure \ref{fig:Fig9} below, 
instead of $t_1$ of Figure \ref{fig:Fig1}. 

\begin{figure}[htbp]
\begin{center}
\includegraphics[width=0.3\linewidth]{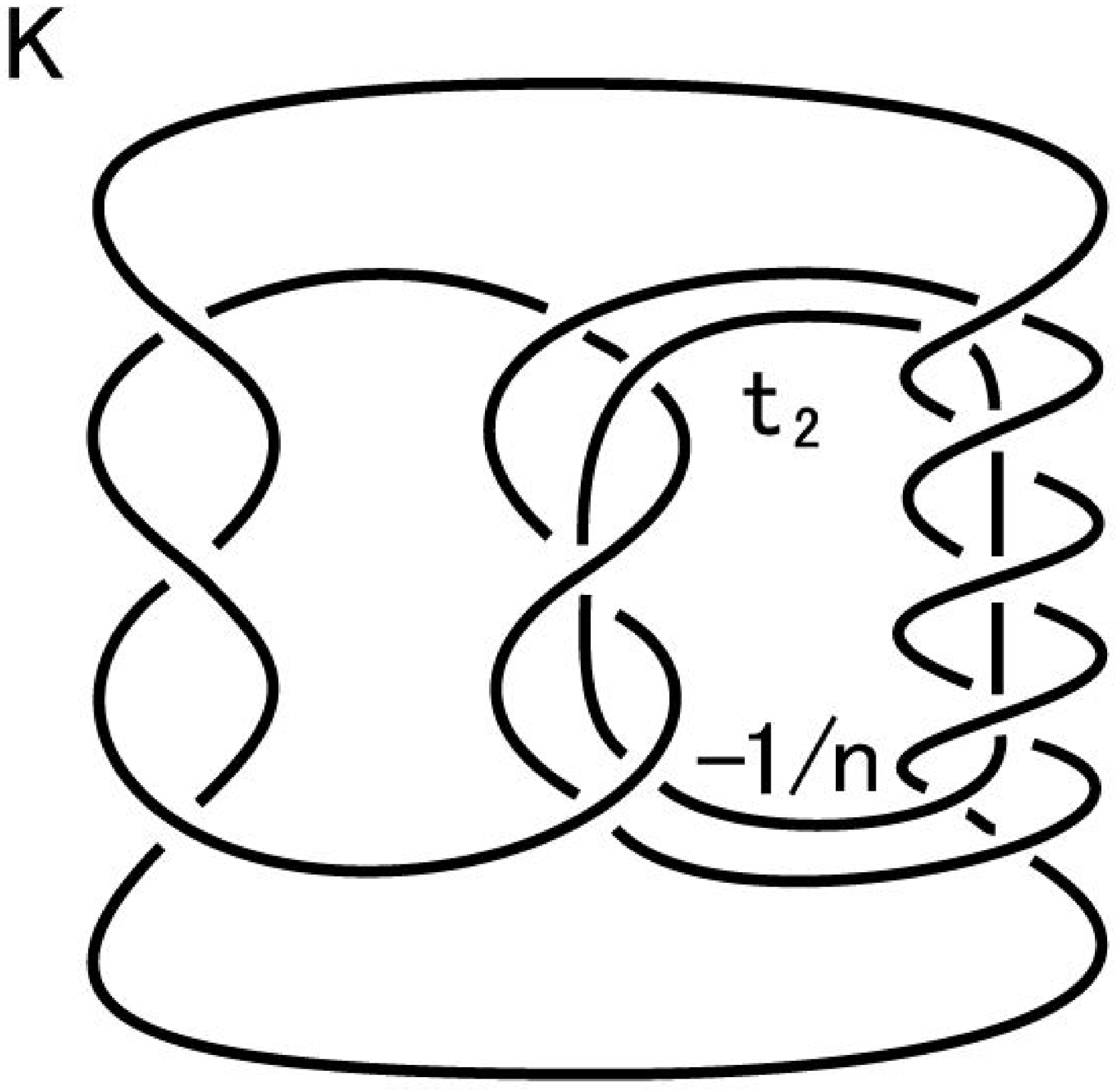}
\caption{}
\label{fig:Fig9}
\end{center}
\end{figure}

Let $K'_n$ be the knot obtained from $K$ 
by doing $n$-twisting along $t_2$. 
Then the argument in the proof of 
Claim \ref{Seifert} shows that $(K'_n ; 1)$ is a 
Seifert fiber space of type $S^2(3, 4, |12n + 5|)$;
see Figures \ref{fig:Fig10}--\ref{fig:Fig13}. 
The arguments in the proofs of 
Claims \ref{Period2}, \ref{hyperbolic}, \ref{non-invertible} and 
\ref{tunnel number} show that the $K'_n$ 
also enjoy the same properties 
as in Example 1,
and that the Seifert fibered surgeries do not come from the 
primitive/Seifert-fibered construction. 

\begin{figure}[htbp]
\begin{center}
\includegraphics[width=1.0\linewidth]{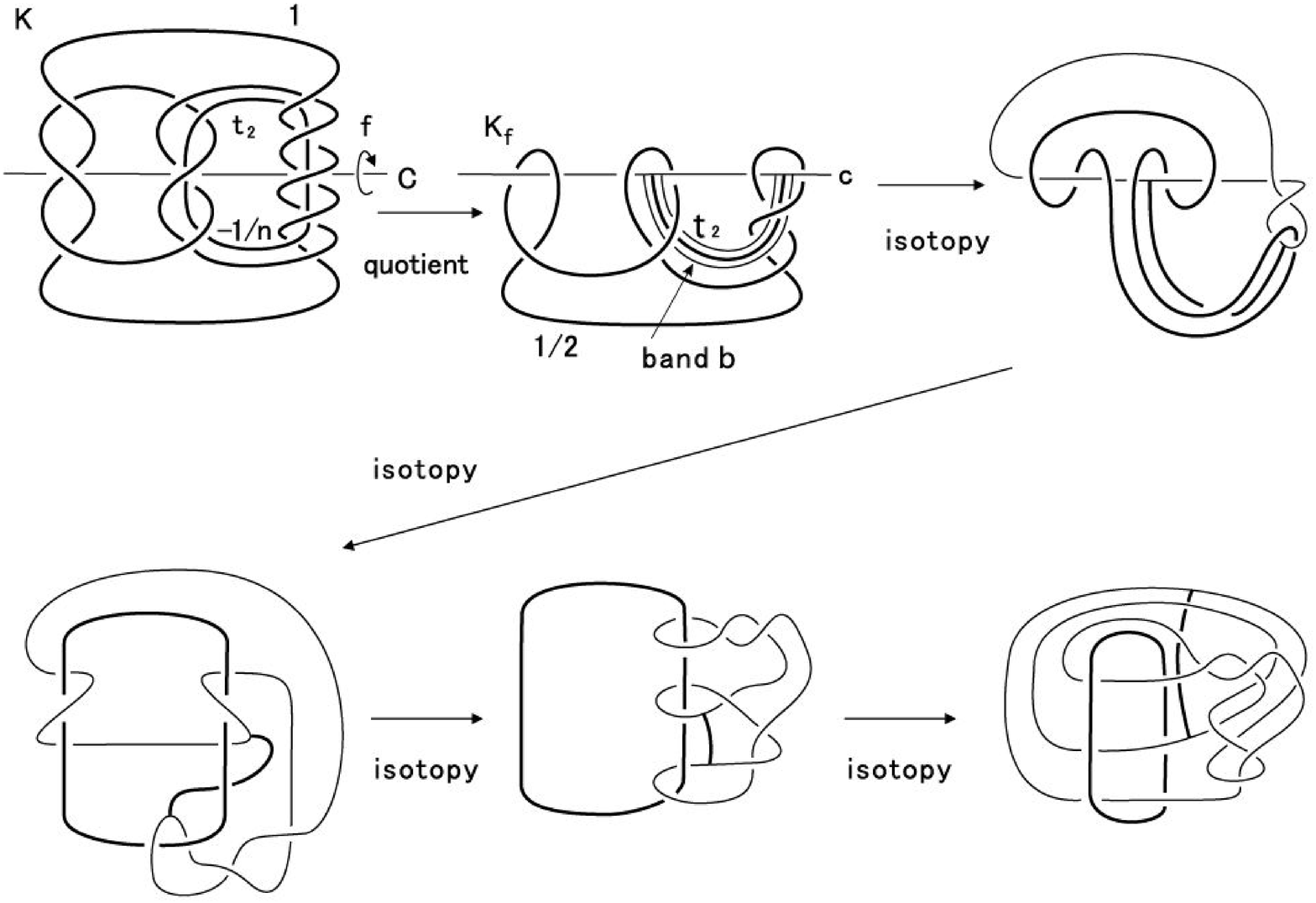}
\caption{}
\label{fig:Fig10}
\end{center}
\end{figure}

\begin{figure}[htbp]
\begin{center}
\includegraphics[width=1.0\linewidth]{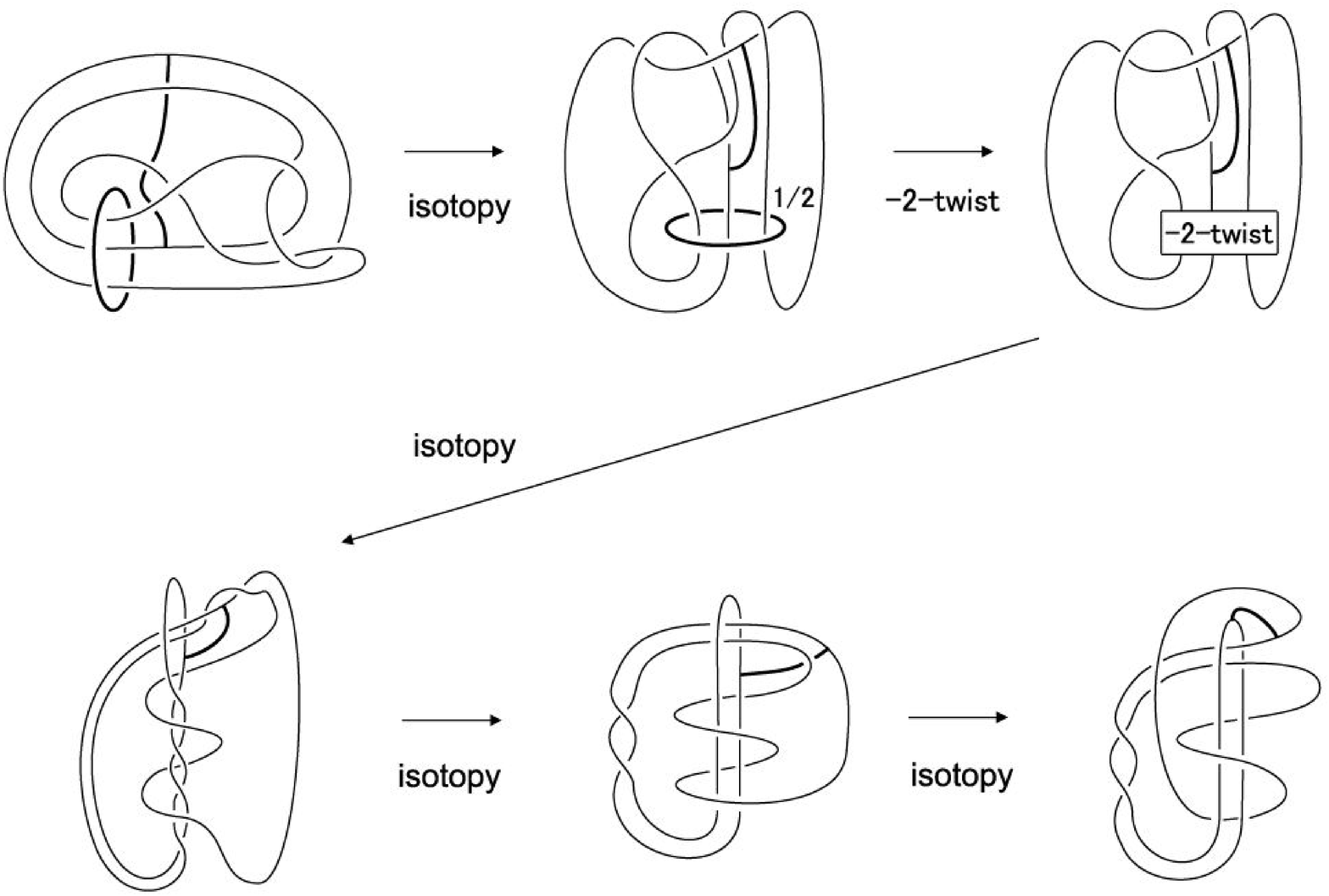}
\caption{Continued from Figure \ref{fig:Fig10}.}
\label{fig:Fig11}
\end{center}
\end{figure}

\begin{figure}[htbp]
\begin{center}
\includegraphics[width=1.0\linewidth]{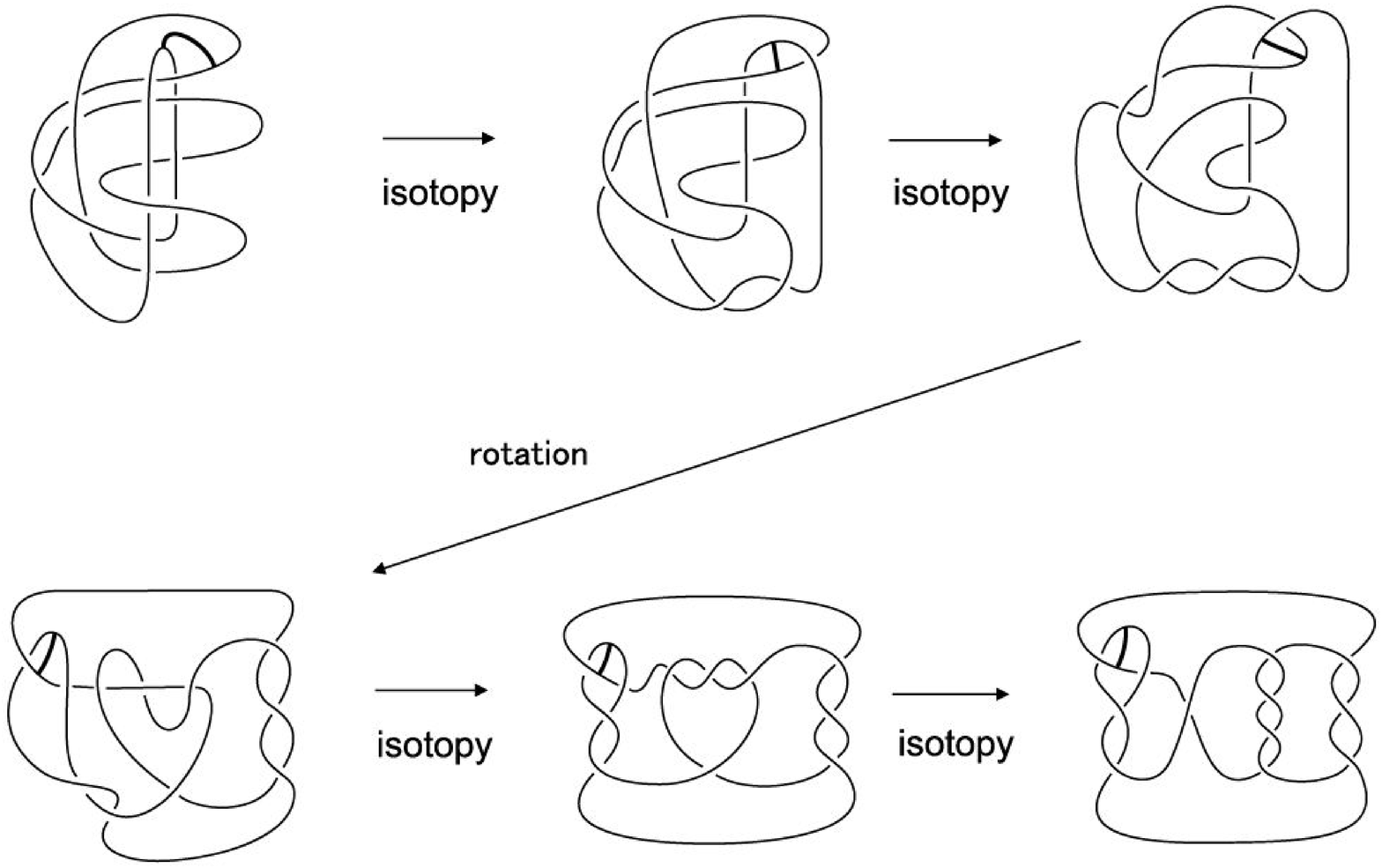}
\caption{Continued from Figure \ref{fig:Fig11}.}
\label{fig:Fig12}
\end{center}
\end{figure}

\begin{figure}[htbp]
\begin{center}
\includegraphics[width=1.0\linewidth]{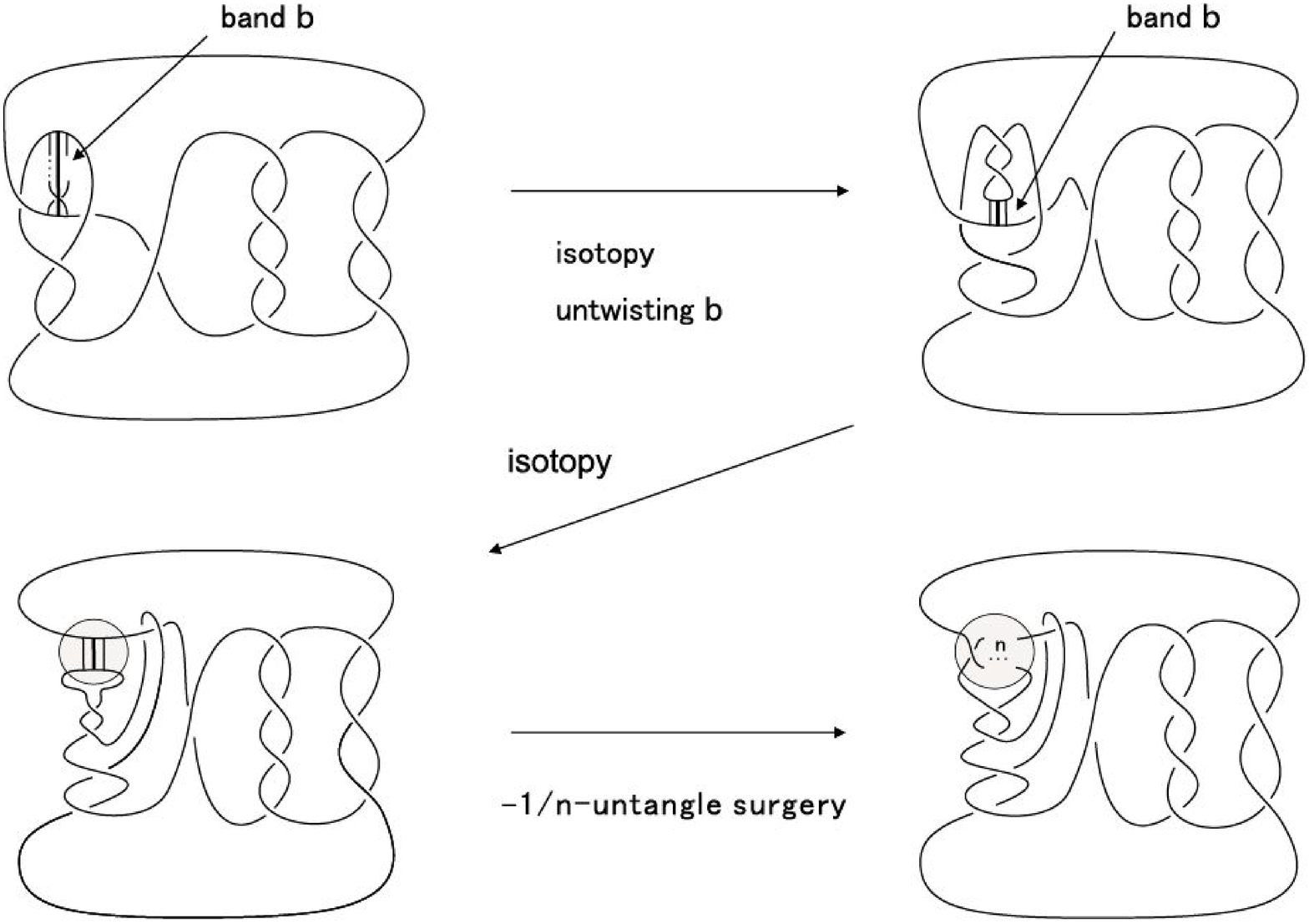}
\caption{Continued from Figure \ref{fig:Fig12}.}
\label{fig:Fig13}
\end{center}
\end{figure}

\section{Remarks and questions}

       In \cite{MM3} it has been conjectured that 
if $(K; r)$ is a Seifert fiber space, 
then it admits a Seifert fibration such that
one of its fibers is unknotted in $($the original$)$ $S^3$. 
For our knots $K_n$ (resp. $K'_n$), 
the trivial knot $t_1^*$ which is 
the dual of $t_1$ (i.e., the core knot of $-1/n$-filling along $t_1$) 
(resp. $t_2^*$ which is the dual of $t_2$)
becomes an exceptional fiber 
of index $|15n + 4|$ in $(K_n ; 1)$ 
(resp. an exceptional fiber of index 
$|12n + 5|$ in $(K'_n ; 1)$). 
Thus the Dehn surgeries described in Examples 1 and 2 
satisfy the conjecture. 
(Song's example mentioned in the 
Introduction also satisfies the conjecture.) \par

        We also mention a geometric aspect of Seifert fibered 
surgeries on hyperbolic knots. 
It was observed in \cite[Section 7]{MM3} that 
short closed geodesics in hyperbolic knot complements 
are often unknotted in $S^3$ and 
become Seifert fibers in the resulting Seifert fiber spaces 
after Dehn surgery. 
An experiment via Weeks' computer program SnapPea \cite{W} 
suggests the table below, 
where $K$ is the $(-3, 3, 5)$-pretzel knot, 
and 
$t_1$, $t_2$ are trivial knots described in Figures 1 and 9. 
Recall that $(K; 1)$ is a Seifert fiber space of type 
$S^2(3, 4, 5)$. 

\medskip

\begin{center}
\begin{tabular}{|c|c|c|c|} \hline
  & $S^3 - K$   & $S^3$ & $(K; 1)$ \\ \hline
$t_1$ & third shortest geodesic & unknot & fiber of index $4$ \\
$t_2$ & shortest geodesic & unknot & fiber of index $5$ \\ \hline
\end{tabular}
\end{center}

\medskip

\noindent
The second shortest geodesic 
is unknotted in $S^3$, 
but it does not become a fiber in 
$(K ; 1)$. In fact it is hyperbolic in $(K; 1)$. 
\par

	We conclude this paper with some questions. 
    Although the knots given in 
Examples 1 and 2 cannot be 
primitive/Seifert-fibered for any genus two Heegaard surface, 
they are still embedded in a genus two Heegaard surface for $S^3$. 
We would like to ask:

\begin{Q}
\label{genus 2}
If $(K;r)$ is a Seifert fiber space, 
then is $K$ embedded in a genus two Heegaard surface for $S^3$? 
\end{Q}

In particular, 

\begin{Q}
\label{tunnel number 2}
If $(K;r)$ is a Seifert fiber space, 
then is the tunnel number of $K$ at most $2$? 
\end{Q}


\bigskip

\begin{thebibliography}{99}
\bibitem{B} J.Berge; 
Some knots with surgeries yielding lens spaces, 
unpublished manuscript. 

\bibitem{Ble} S. A. Bleiler; 
Knots prime on many strings, 
Trans.\ Amer.\ Math.\ Soc.\ \textbf{282} (1984), 385--401. 

\bibitem{Bl} S. A. Bleiler;
Prime tangles and composite knots, 
Lect.\ Notes in Math.\ vol.\ \textbf{1144}, Springer-Verlag, 1985, 
pp. 1--13. 

\bibitem{BH} S. Bleiler and C. Hodgson;
Spherical space forms and Dehn filling,
Topology  \textbf{35} (1996), 809--833.

\bibitem{BZ} S. Boyer and X. Zhang;
Finite surgery on knots,
J. Amer.\ Math.\ Soc.\ \textbf{9} (1996), 1005--1050.

\bibitem{BuZ} G. Burde and H. Zieschang; 
Knots, de Gruyter Studies in Mathematics \textbf{5}, 1985.

\bibitem{D} J. Dean;
Hyperbolic knots with small Seifert-fibered Dehn surgeries,
Ph.D. thesis, University of Texas at Austin, 1996.

\bibitem{D2} J. Dean; 
Small Seifert-fibered Dehn surgery on hyperbolic knots, 
Algebraic and Geometric Topology \textbf{3} (2003), 
435--472. 

\bibitem{EM1} M. Eudave-Mu\~noz;
Non-hyperbolic manifolds obtained by Dehn surgery on a hyperbolic knot,
In: Studies in Advanced Mathematics vol.\ \textbf{2}, part~1,
(ed.\ W. Kazez), 1997, Amer.\ Math.\ Soc.\ and International Press,
pp. 35--61.

\bibitem{EM2} M. Eudave-Mu\~noz; 
On hyperbolic knots with Seifert fibered Dehn surgeries, 
Topology Appl.\ \textbf{121} (2002), 
119--141. 

\bibitem{Ga} D. Gabai; Foliations and the topology of 3-manifolds
{\rm III},
J.\ Diff.\ Geom.\ \textbf{26} (1987), 479--536.
\bibitem{GS} F. Gonz\'alez-Acu\~na and H. Short; 
Knot surgery and primeness, 
Math.\ Proc.\ Camb.\ Phil.\ Soc.\ \textbf{99} (1986), 
89--102.

\bibitem{Go} C. McA. Gordon; Dehn Filling; a survey, 
Knot theory (Warsaw, 1995), 129--144,
Banach Center, Publ. 42, Polish Acad. Sci., Warsaw, 1998.

\bibitem{GL} C. McA. Gordon and J. Luecke;
Knots are determined by their complements,
J. Amer.\ Math.\ Soc.\ \textbf{2} (1989), 371--415.

\bibitem{Mat} T. Mattman; 
The Culler-Shalen seminorms of pretzel knots, 
Ph.D. thesis, 
McGill University, Montr\'eal, 2000.

\bibitem{MM2} K. Miyazaki and K. Motegi;
Seifert fibered manifolds and Dehn surgery {\rm II},
Math. Ann. \textbf{311} (1998), 647--664.

\bibitem{MM3} K. Miyazaki and K. Motegi;
Seifert fibered manifolds and Dehn surgery {\rm III},
Comm.\ Anal.\ Geom.\ \textbf{7} (1999), 551--582.

\bibitem{MM7} K. Miyazaki and K. Motegi;
On primitive/Seifert-fibered constructions,
Math.\ Proc.\ Camb.\ Phil.\ Soc.\ (to appear). 

\bibitem{Mon} J. M. Montesinos; 
Surgery on links and double branched coverings of $S^3$, 
Ann.\ Math.\ Studies \textbf{84} (1975), 227--260.

\bibitem{Mor} K. Morimoto; 
On the additivity of h-genus of knots, 
Osaka J. Math. \textbf{31} (1994), 137--145.

\bibitem{Mor2} K. Morimoto;
There are knots whose tunnel numbers go down under connected sum,
Proc. Amer. Math. Soc. \textbf{123} (1995), 3527--3532. 

\bibitem{Mot} K. Motegi; Dehn surgeries, group actions and 
Seifert fiber spaces, 
Comm.\ Anal.\ Geom.\ \textbf{11} (2003), 343--389. 

\bibitem{S} H. Schubert; Knoten und Vollringe, 
Acta Math.\ \textbf{90} (1953), 131--286.

\bibitem{W} J. Weeks; SnapPea: 
a computer program for creating and studying hyperbolic $3$-manifolds, 
freely available from \hfil\break
http://thames.northnet.org/weeks/index/SnapPea.html

\bibitem{Z} H. Zieschang; On simple systems of paths on complete 
pretzels, 
Amer.\ Math.\ Soc.\ Transl.\ \textbf{92}, 127--137.

\end{thebibliography}
\end{document}